\documentclass{article}
\usepackage{latexsym}
\usepackage{amsthm}  % for proof environment
\usepackage[all]{xy}
\usepackage{stmaryrd}
\usepackage{amssymb}
\usepackage{amsmath}
\usepackage{mathptmx}

\newcommand{\mbb}[1]{\mathbf{#1}}

\newcommand{\msc}[1]{\mathcal{#1}}
\newcommand{\mrm}[1]{\mathrm{#1}}

\newcommand{\C}{\mbb{C}}

\newcommand{\Z}{\mbb{Z}}

\newcommand{\Q}{\mbb{Q}}

\newcommand{\A}{\mbb{A}}

\newcommand{\Pro}{\mbb{P}}

\newcommand{\tens}{\varotimes}

\newcommand{\Aut}{\mrm{Aut~}}

\newcommand{\ses}[5]{{\xymatrix{0\ar[r]&{#1}\ar[r]^{#2}&{#3}\ar[r]^{#4}&{#5}\ar[r]& 0}}}

\newcommand{\sh}[1]{{\msc{#1}}}

\newcommand{\shext}{{\msc{E}\!\!\mathit{xt}}}

\newcommand{\Ext}{\mrm{Ext}}

\newcommand{\Spec}{\mrm{Spec~}}

\newcommand{\Def}{\mrm{Def}}

\theoremstyle{definition}
\newtheorem{defn}{Definition}[section]
\newtheorem{rmk}[defn]{Remark}

\newtheorem{exa}[defn]{Example}
\theoremstyle{plain}
\newtheorem{thm}[defn]{Theorem}

\newtheorem{lem}[defn]{Lemma}
\newtheorem{prop}[defn]{Proposition}

\date{}		% Fill in date or delete command for todayÕs date.

\begin{document}

\title{Stable degenerations of surfaces isogenous to a product of curves}
\author{Michael A. van Opstall}
\maketitle

\section{Introduction}

In \cite{cat:fs} and its sequel, \cite{cat:corr}, Catanese studies surfaces
which admit an unramified cover by a product of curves of genus greater than
one. Then main interesting result from the point of view of moduli theory is
that the moduli space of such surfaces (that is, of canonically polarized
surfaces with the same numerical invariants $K^2$ and $\chi(\sh{O})$) is 
either irreducible or has two connected components, swapped by complex
conjugation.

In \cite{ksb:3f}, Koll\'ar and Shepherd-Barron introduced a compactification
of the moduli space of canonically polarized surfaces. The details of the
construction were later filled in by several authors. In particular, they
give a {\em stable reduction} procedure, by which any one-parameter family
of surfaces over a punctured disk can be completed to a family of so-called
{\em stable surfaces} over a finite cover of the disk. The stable reduction
is obtained by taking the relative canonical model of a semistable
resolution of the original family. 

This article is concerned with determining the stable surfaces which occur
at the boundary of this moduli space. The moduli space of stable curves 
provides us with candidates for the boundary surfaces, and we simply verify
that these candidates are already stable, so the full force of the 
minimal model program is not required. In contrast, a forthcoming article
demonstrates the use of flips in determining the stable degenerations of
symmetric squares of curves.

Pairing our results with Catanese's results, we obtain that the two components
of the moduli space of canonically polarized surfaces do not meet in the
stable compactification. To my knowledge, this is the first example of 
disconnectedness of a moduli space of {\em stable} surfaces after fixing the
diffeomorphism class of a smooth member. The disconnectedness of the moduli
space of canonically polarized surfaces with fixed differentiable structure
was first proved by Manetti. However, to establish that his surfaces which
lie on different components of the moduli space are diffeomorphic, he 
degenerates them all to a single stable surface. Catanese's examples, on the
other hand, are obviously diffeomorphic surfaces.

\section{Stable surfaces and surfaces isogenous to a product}

First we define the higher-dimensional analogue of the nodes which are 
allowed on stable curves:

\begin{defn}
A surface $S$ has {\em semi-log canonical} (slc) singularities if
\begin{enumerate}
\item $S$ is Cohen-Macaulay;
\item $S$ has normal crossings singularities in codimension one;
\item $S$ is $\Q$-Gorenstein, i.e. some reflexive power of the dualizing
sheaf of $S$ is a line bundle;
\item for any birational morphism $\pi:X\rightarrow S$ from a smooth variety,
if we write (numerically)
\[
K_X=\pi^*K_S+\sum a_iE_i,
\]
then all the $a_i\geq 1$.
\end{enumerate}
\end{defn}

The first condition ensures that the dualizing sheaf exists; the second
gives that it is an invertible sheaf off a subset of codimension 2, so it
can be extended to give a Weil divisor class $K_S$. The third condition states
that some multiple of this class is Cartier, so we can make sense of the
formula occuring in the fourth condition. A complete classification of
slc surface singularities can be found in \cite{ksb:3f}. Suppose $S$ is
a $\Q$-Gorenstein surface. Let $S^\nu$ be the normalization, and $D$ the 
inverse image of the codimension 1 singular set under the normalization
morphism. Then the condition that $S$ be slc is equivalent to the condition
that $(S^\nu,D)$ is a log canonical (lc) pair.

\begin{defn}
A {\em stable surface} is a projective, reduced surface $S$ with slc 
singularities such that some reflexive power $\omega_S^{[N]}$ is an ample
line bundle. The smallest such $N$ such that $\omega_S^{[N]}$ is a line
bundle is called the {\em index} of $S$. A family of stable surfaces
is a flat morphism $X\rightarrow B$ whose fibers are stable surfaces and
whose relative dualizing sheaf $\omega_{X/B}$ is $\Q$-Cartier.
\end{defn}

After fixing a Hilbert polynomial $P$, there is a bound for the index of
a stable surface with Hilbert polynomial $P$. This is one of the ingredients
in the construction of the moduli space of stable surfaces with fixed
Hilbert polynomial $P$, which compactifies the moduli space of canonically
polarized surfaces (after possibly throwing away some components parameterizing
only singular surfaces). This moduli space is proper and separated. 
The moduli space would not be separated without the requirement that families
of stable surfaces be relatively $\Q$-Gorenstein. It is worth noting that
in a later article (\cite{k:pcm}) Koll\'ar strengthened the conditions that
a family of stable surfaces should satisfy. It is unknown whether the 
$\Q$-Gorenstein assumption alone implies these stronger condition. For our 
purposes, the weaker condition is sufficient, since it is known to imply the 
stronger conditions in the case of one-parameter families whose general
member is smooth.

We now recall some of the definitions from \cite{cat:fs}.

\begin{defn}
A surface $S$ is called isogenous to a product if it admits an unramified
cover by a product of curves of genus two or higher.
\end{defn}

\begin{rmk}
Catanese calls such a surface isogenous to a {\em higher} product, but we will
not be interested in surfaces covered by, say, a product of elliptic curves.
It is clear that any such surface is canonically polarized, since the cover
by a product of curves contains no rational curves, so $S$ contains none.
\end{rmk}

The following is a summary of 3.10-3.13 of {\em loc. cit.}:

\begin{prop}
A surface $S$ isogenous to a product can be written uniquely as 
$(C_1\times C_2)/G$ for some group $G$ which embeds into $\Aut C_1$ and
$\Aut C_2$, as long as $C_1$ and $C_2$ are not isomorphic.
If $C_1$ and $C_2$ are isomorphic, the subgroup of $G$ consisting of 
automorphisms not switching the factors is required to embed into the
automorphism group of each factor to obtain this {\em minimal realization}.
\end{prop}

This proposition allows us to describe small deformations of surfaces 
isogenous to a product. This was done in {\em loc. cit.}, but we present
a simple and purely algebraic proof here. 

\begin{defn}
The functor of deformations of a stable curve $C$ with the action of a finite
group $G$ assigns to an artin ring $A$:
\begin{enumerate}
\item a flat morphism $X\rightarrow \Spec A$,
\item an embedding of $G$ in the group $\Aut_{\Spec A} X$ of automorphisms
of the family over the base,
\item and an equivariant isomorphism of the special fiber of the family
$X$ with the stable curve $C$.
\end{enumerate}
\end{defn}

\begin{prop}
The functor of deformations of a stable curve $C$ together with a subgroup 
$G$ of the automorphism group is ``well-behaved'' (i.e. satisfies the 
Schlessinger conditions) and unobstructed, with tangent space 
$\Ext^1(\Omega_C,\sh{O}_C)^G$. Such pairs have a proper moduli space,
finite over a closed subvariety of the moduli space of stable curves. 
\end{prop}

The proof of this may be found in \cite{gmg1}. Note that the statements here 
depend on the exactness of taking invariants by a finite group in 
characteristic zero, and that the moduli space mentioned above is not
proper in positive characteristic because of wild ramification.

\begin{prop}\label{locmod}
The Kuranishi space of a surface $S$ minimally realized as a free quotient
$(C_1\times C_2)/G$ is (locally analytically) isomorphic to the product
of Kuranishi spaces of the pairs $(C_1,G)$ and $(C_2,G)$. Consequently, the
Kuranishi space is smooth, and the moduli space of canonically polarized 
surfaces is irreducible at $S$.
\end{prop}

\begin{proof}
Since all deformation functors in question are ``nice'' enough, we may 
resolve the question by checking to first order:
\begin{eqnarray*}
H^1(S,\sh{T}_S)&=& H^1(C_1\times C_2,\sh{T}_{C_1\times C_2})^G \\
&=& H^1(C_1,\sh{T}_{C_1})^G\oplus H^1(C_2,\sh{T}_{C_2})^G.
\end{eqnarray*}
The first line is valid only when $G$ acts freely, but the second line
works generally for products of canonically polarized varieties, assuming
$G$ acts on both factors (see e.g., \cite{vo:mpc}).

Since deformations of curves with group action are unobstructed, the
Kuranishi space of $S$ is smooth. Since $S$ has a finite automorphism
group, the moduli space of $S$ locally near $S$ has only finite
quotient singularities, so cannot be reducible.
\end{proof}

\section{Degenerations}

Our main goal is the following theorem:

\begin{thm}
Suppose $X\rightarrow \Delta'$ is a family of surfaces isogenous to products
over a punctured disk.
Then possibly after a finite change of base, totally ramified over the origin
in the disc, $X$ (or a pullback thereof) can be completed to a family of
stable surfaces over the disk whose central fiber is a quotient of a product
of stable curves (possibly by a non-free group action).
\end{thm}

\begin{proof}
By \ref{locmod}, we may assume that $X$ is of the form 
$(Y_1\times_{\Delta'} Y_2)/G$ where $Y_1$ and $Y_2$ are families of smooth
curves with $G$-action, such that the $G$-action is fiberwise and free on 
$Y_1\times_{\Delta'} Y_2$. Since the moduli functor of stable curves with 
automorphism group $G$ is proper, after a base change (which we will suppress
in our notation), we obtain a family $\tilde{X}$ of the desired form.

It remains to see that the central fiber is a stable surface, and that the
family $\tilde{X}$ is a family of stable surfaces. $\tilde{X}$ is obtained
by taking the quotient of a family $\tilde{Y}$ of stable surfaces by a group
action. Since the group acts freely on the general fiber, the quotient
morphism is \'etale in codimension 1. In this case, \cite{km:bgav}, 
Proposition 5.20 ensures that $\tilde{X}$ is $\Q$-Gorenstein, so the special
fiber is as well. Well known results ensure that the special fiber is
Cohen-Macaulay. A finite quotient of a variety which is normal crossings
in codimension one is normal crossings in codimension one by Corollary 1.7 
of \cite{aal:wncs}. Now \cite{km:bgav} 5.20 (appropriately modified to
take into account non-normal varieties) states that a finite quotient of an slc
variety is slc as soon as it is $\Q$-Gorenstein an normal crossings in
codimension 1.

Finally, since the relative canonical sheaf $\omega_{\tilde{Y}/\Delta}$ is
ample since $Y$ is a family of products of curves of general type, and
the quotient morphism is unramified in codimension one, the relative
canonical sheaf of $\tilde{X}$ is also ample, so $\tilde{X}$ is a family
of stable surfaces, hence the special member is uniquely determined by
separatedness of the moduli space of stable surfaces.
\end{proof}

\section{Digression: Deformation of curves with group action}

For a family $X\rightarrow S$ of proper varieties, the functions
$\dim~T^n(X_s)$ are Zariski upper semicontinuous; that is, they may
jump up at certain points. This is proved in \cite{pal:dcs}, where actually
a stronger result is obtained:

\begin{thm}[Palamodov]
Let $X\rightarrow S$ be a flat proper morphism.
The functions $\sum_{i=0}^{2m}(-1)^i\dim~T^{n+i}(X_s)$ are upper 
semicontinuous, for $m>0$ and $n\geq 0$.
\end{thm}

A consequence for stable curves is that $\dim~T^1(C_s)$ is constant in a 
family of stable curves, since $T^0$ and $T^2$ vanish for stable curves. This 
is one way of showing that even for a stable curve, $T^1(C)$ is 
$3g-3$-dimensional, although direct computation is straightforward, even
for an arbitrary nodal curve.

Does such a strong semicontinuity property hold for the functors $T^i(-,G)$?
A weak version follows from Palamodov's theorem, at least in our situation of
finite groups and varieties defined over fields of characteristic zero:

\begin{lem}
Let $X\rightarrow S$ be a family of varieties with $G$-action. Then 
$\dim~T^i(X_s,G)$ is an upper semicontinuous function of $s\in S$. 
\end{lem}

\begin{proof}
From \cite{pal:dcs}, we have upper semicontinuity for 
$\dim~T^i(X_s)$. Under our assumptions, $T^i(X_s,G)=T^i(X_s)^G$. The set of
$G$-invariants is the intersection of the kernels of the linear maps
$g-I$ for all $g\in G$, and it is well known that the dimension of the
kernel of a vector bundle map is upper semicontinuous.
\end{proof}

A review of the computation of of $T^1(C)$ for a stable curve will be useful
in computing $T^1(C,G)$. First of all, for nodal curves (or more generally,
complete intersection curves), the local-to-global Ext spectral sequence
computes $T^i$: $T^i(C)=0$ for $i>1$, $T^0(C)=H^0(C,\sh{T}_C)$, and we can
write, using the exact sequence of low degree terms of this spectral sequence,
\[
\ses{H^1(\sh{T}_C)}{}{\Ext^1(\Omega^1_C,\sh{O}_C)}{}{H^0(C,\shext^1(\Omega^1_C,
\sh{O}_C))}
\]
an extension describing $T^1$. Let us describe the topological data of
a stable curve as in \cite{hm:mc}: let $\nu$ be the number of irreducible
components. Then let $g_i$ be the genus of the normalization of the $i$th
irreducible component for $i=1,\ldots,\nu$. Finally, let $\delta$ equal the
number of nodes. Then the (arithmetic) genus of $C$ is given by
\[
g=\sum_{i=1}^\nu g_i+\delta-\nu+1.
\] 

The space $H^0(C,\shext^1(\Omega^1_C,\sh{O}_C))$ is easy to compute: the 
sheaf $\shext^1(\Omega^1_C,\sh{O}_C)$ is supported at 
the nodes with one-dimensional stalks. Therefore, the dimension of
$H^0(C,\shext^1(\Omega^1_C,\sh{O}_C))$ is $\delta$.

For a nodal curve $C$, let $C^\nu$ be the normalization. Then there are
pairs $(p_k,q_k)$ of points for $k=1,\ldots,\delta$ lying over each node
and giving an exact sequence of sheaves on $C^\nu$:
\[
\ses{\sh{T}_{C^\nu}(-\sum_{k=1}^\delta(p_k+q_k))}{}{\sh{T}_{C^\nu}}{}
{\sh{T}_{C^\nu}\tens\sh{O}_{\sum_{k=1}^\delta(p_k+q_k)}}
\]

By stability, $H^0(\sh{T}_{C^\nu}(-\sum_{k=1}^\delta(p_k+q_k)))=0$, and
for dimension reasons, $H^1(\sh{T}_{C^\nu}\tens
\sh{O}_{\sum_{k=1}^\delta(p_k+q_k)})=0$.Therefore 
\[
\dim~T^1(C)=\delta+h^0(\sh{T}_{C^\nu}\tens\sh{O}_{\sum_{k=1}^\delta(p_k+q_k)})
-\chi(\sh{T}_{C^\nu}),
\]
where we use the fact that $H^1(C,\sh{T}_C)=H^1(C^\nu,
\sh{T}_{C^\nu}(-\sum_{k=1}^\delta(p_k+q_k)))$.

The pieces are easy to compute:
\[
-\chi(\sh{T}_{C^\nu})=3\sum_{i=1}^\nu g_i-3\nu
\]
and 
\[
h^0(\sh{T}_{C^\nu}\tens\sh{O}_{\sum_{k=1}^\delta(p_k+q_k)})=2\delta.
\]
We conclude that
\[
\dim~T^1(C)=3\sum_{i=1}^\nu g_i+3\delta-3\nu
\]
which is $3g-3$ in light of the formula given above for $g$.

Now to the equivariant case: since the operation of taking $G$-invariants
is exact, we obtain a similar formula:
\[
\dim~T^1(C,G)=h^0(\shext^1(\Omega^1_C,\sh{O}_C))^G
+h^0(\sh{T}_{C^\nu}\tens\sh{O}_{\sum_{k=1}^\delta(p_k+q_k)})^G
-\chi(\sh{T}_{C^\nu})^G
\]
where by $\chi(\sh{T}_{C^\nu})^G=\dim~H^1(\sh{T}_{C^\nu})^G-
\dim~H^0(\sh{T}_{C^\nu})^G$. Note that the $G$-action extends to 
$C^\nu$, since $G$ acts birationally on $C^\nu$ (normalization is
birational) and $C^\nu$ is a smooth curve, so any birational morphism
extends to an automorphism. 

The following example illustrates the method of computation in the case
where $G$ does not permute the components of $C^\nu$. In this case,
$C^\nu/G$ is a disjoint union of smooth curves. We also make use of the fact
that the deformation functor $\Def_{(C,G)}$ of a smooth curve with group
action is isomorphic to the functor $\Def_{(C/G,B)}$ where $B$ is the 
branch locus of the cover $C\rightarrow C/G$.

\begin{exa}
Let $C$ be the curve defined by the homogenous equation $xyz^2+x^4+y^4=0$
in $\Pro^2$. The arithmetic genus of $C$ is 3, and $C$ has a single node
at the point with homogeneous coordinates $(0:0:1)$. The surface
$xyz^2+x^4+y^4+tz^4$ in $\Pro^2_{\A_t^1}$ is a one-parameter smoothing.

Let $G=\Z_2$ act on $C$ by swapping $x$ and $y$ and sending $z$ to $-z$.
The fixed locus of this automorphism of $\Pro^2$ is the set of points
with homogeneous coordinates $(1:-1:z)$ and the points $(0:0:1)$ and 
$(1:1:0)$, that is, a line together with an isolated point. Therefore, 
$G$ acts on a general member of the smoothing (a smooth quadric) fixing
four points. The action of $G$ on $C$ fixes three points, one of them the 
node. 

Let $C^\nu$ denote the normalization of $C$ and $p$ and $q$ be the
two points of $C^\nu$ lying over the node of $C$. The group action fixes
the node of $C$, but swaps the two branches. This means that the action
lifted to $C^\nu$ swaps $p$ and $q$. Analytically locally near the node,
the curve looks like $\C[x,y]/(xy)$ together with the action of $G$ swapping
$x$ and $y$. Since the $G$ action lifts to a versal deformation, we 
conclude that every element of $H^0(\shext^1(\Omega^1_C,\sh{O}_C))$ is
$G$-invariant. This vector space is one-dimensional.

Next we compute $H^0(\sh{O}_{p,q}\tens\sh{T}_{C^\nu})^G$. The space
$H^0(\sh{O}_{p,q}\tens\sh{T}_{C^\nu})$ is a two-dimensional vector space.
Since the $G$-action on $C^\nu$ swaps $p$ and $q$, only the ``diagonal'' 
(in an appropriate basis) is fixed by $G$. So this space contributes one
dimension to $T^1(C,G)$.

Finally, we need to compute $\chi(\sh{T}_{C^\nu})^G$. $C^\nu$ is a genus
2 curve (genus drops by the number of nodes upon normalization) and it
double covers $C^\nu/G$ with two points fixed (the two points of $C$ other
than the node which are fixed by $G$). Therefore, $C^\nu/G$ is an
elliptic curve with two marked points. It follows that
$\chi(\sh{T}_{C^\nu})^G=-2$, and therefore $T^1(C,G)$ is four-dimensional.

In this example, there is no jumping of dimension: the general fiber
is a genus 3 curve with $G$-action whose quotient is a 4-pointed
elliptic curve. Therefore $T^1$ of the generic fiber with $G$-action
is also four-dimensional.
\end{exa}

Note that it is critical how $G$ acts near a node in order to
compute $H^0(\shext^1(\Omega^1_C,\sh{O}_C))^G$; it is not the case
that if $G$ fixes a node, then every element of 
$H^0(\shext^1(\Omega^1_C,\sh{O}_C))$ is fixed.

Now we show that the example is general: $\dim~T^1(C_s,G)$ cannot jump
in a family of stable curves.

\begin{thm}
Let $C\rightarrow S$ be a family of stable curves with an action of a
finite group $G$. Then $\dim~T^1(C_s,G)$ is locally constant.
\end{thm}

\begin{proof}
First of all, the deformation functor of curves of genus greater than one 
with group action is formally smooth (this is true even in positive
characteristic, assuming that the action is tame, by \cite{gmg1}, 5.1). This
functor has a versal deformation, which is universal, since $T^0(C,G)=0$
for stable curves. The universal family induces a morphism to the moduli
space, and the induced morphism is finite and surjective onto a neighborhood
of the point of the moduli space corresponding to the isomorphism class
of $(C_0,G)$. Therefore the moduli space has at worst finite quotient
singularities, and its dimension is equal to $\dim~T^1(C_0,G)$. Since the
moduli space is locally irreducible, its dimension cannot jump at special
points, from which the theorem follows.
\end{proof}

\section{Application to Catanese's examples}

In \cite{cat:corr}, Catanese gives a family of examples of moduli spaces
of smooth surfaces with fixed $K^2$ and $\chi$ fixed which have two components
interchanged by complex conjugation. We review his construction here and
study the degenerations of his surfaces.

The construction of the example begins with the construction of a
triangle curve (i.e. a Galois cover of $\Pro^1$ branched at 3 points)
which is not antiholomorphic to itself. Let $C$ denote this curve and
$G$ denote the Galois group of the cover $C\rightarrow \Pro^1$. $G$ is
therefore a quotient of the fundamental group of $\Pro^1$ minus three points,
and is consequently generated by 2 elements. Choose $h\geq 2$ and a curve
$C_1'$ of genus $h$. Then the fundamental group of $C_1'$ surjects onto
$G$, so there exists an \'etale cover $C_1\rightarrow C_1'$ with Galois
group $G$. Then the surface $S=(C_1\times C)/G$ is isogenous to a product
of curves of general type (the triangle curve constructed is not the
elliptic curve with $j$-invariant 1728, which is the only triangle curve
not of general type). 

The critical result for finding multiple components of the moduli space is 
Catanese's Proposition 3.2: the existence of an antiholomorphic isomorphism of 
two surfaces minimally realized as surfaces isogenous to products of curves of 
general type implies antiholomorphic isomorphisms of the factors (up to
reordering the factors). In what follows, denote by $\overline{X}$ the 
complex conjugate of the manifold $X$. 

Choosing any $C_2'$ of genus $h$ and a cover $C_2$ of $C_2'$ with Galois
group $G$ as above, suppose $(C_1\times C)/G\cong \overline{(C_2\times C)/G}$.
Then there is an antiholomorphic isomorphism of $(C_1\times C)/G$ with
$(C_2\times C)/G$, and hence, an antiholomorphic automorphism of $C$, which
is impossible by the construction of $C$. Catanese claims that the various
choices of $C_2$ fill out a component of the moduli space. But 
$\overline{(C_2\times C)/G}$ is diffeomorphic to $(C_2\times C)/G$, and hence
also has a point in the moduli space, which cannot be on this component.
Therefore the moduli space has at least two components.

Now let us consider the stable degenerations of these surfaces, and address
the question of whether the two components are joined together by deformations
through stable surfaces.
The results in this chapter show that the (small) deformations of $S$ are just
the $G$-equivariant deformations of $C_1$, or equivalently, the deformations
of $C_1/G$. Let $\overline{M}$ denote the moduli space of smoothable stable
surfaces occuring as degenerations of $(C_1\times C)/G$ or its conjugate.
$\overline{M}$ has two irreducible components; is it connected?

Suppose $\overline{M}$ were connected: then there would exist a surface
$(C'\times C)/G$ on the boundary of the moduli space which lies on the
closure of both components. Since both components come from curves with 
$G$-action, the map induced from the Kuranishi space of $(C'\times C,G)$
must surject onto a neighborhood of the corresponding boundary point. But 
the Kuranishi space of $(C'\times C, G)$ is irreducible, so it cannot map
onto two components. So the disconnection of various moduli spaces considered
in \cite{cat:corr} continues in the stable compactification.

Note that this argument is not strong enough in general to claim that a 
moduli space of surfaces isogenous to a product of curves is always irreducible
at the boundary: it just rules out deformations to other surfaces isogenous
to a product with the same Galois group. By the results of \cite{cat:corr},
there are at most two components of the moduli space (after fixing Hilbert
polynomial) parameterizing smooth varieties, but there may be a component
parameterizing only singular surfaces meeting both other components along
the boundary.

\end{document}